\documentclass[11pt]{article}
\usepackage{amsmath,amsfonts,epsfig}
\usepackage{mathptmx}

\usepackage{amsmath,amsfonts,latexsym,amssymb}
\usepackage{mathrsfs}
\usepackage{verbatim}
\usepackage[utf8]{inputenc}
\usepackage[T1]{fontenc}
\usepackage{ae,aecompl}
\usepackage{braket}
\usepackage{color}
\usepackage{euscript}

\renewcommand{\bold}[1]{\medskip \noindent {\bf \boldmath #1
                        }\nopagebreak[4]}

\renewcommand{\Im}{\operatorname{Im}}

\newtheorem{theorem}{Theorem}[section]

\newtheorem{lemma}[theorem]{Lemma}

\newtheorem{corollary}[theorem]{Corollary}

\renewcommand{\hbar}{\bar{{\mathbb H}}^3}

\newcommand{\CC}{\mathbb C}

\newcommand{\R}{\mathbb R}
\newcommand{\Q}{\mathbb Q}
\newcommand{\N}{\mathbb N}
\newcommand{\Z}{\mathbb Z}

\newcommand{\Hp}{{\mathbb H}^2}
\newcommand{\Hs}{{\mathbb H}^3}

\def\Rog{\mathcal L}
\def\eproof{$\Box$ \medskip}

\renewcommand\marginpar[1]{} 

\begin{document}

\title{\bf \Large Dilogarithm identities for solutions to Pell's equation in terms of continued fraction convergents } \author{Martin
   Bridgeman\thanks{M. Bridgeman's research was supported by NSF grants DMS-1500545, DMS-1564410.} }

\date{\today}

\maketitle
\begin{abstract}
In this paper we give describe a new connection between the dilogarithm function and solutions to Pell's equation $x^2-ny^2 = \pm 1$. For each solution $x,y$ to Pell's equation we obtain a dilogarithm identity whose terms are given by the continued fraction expansion of the associated unit $x+y\sqrt{n} \in \Z[\sqrt{n}]$.  We further show that Ramanujan's dilogarithm value-identities correspond to an identity for the regular ideal hyperbolic hexagon.\end{abstract}

\section{Dilogarithm and Pell's Equation}
{\bf Dilogarithm:} The dilogarithm function $Li_2(z)$ is given by the power series
$$Li_2(z) = \sum_{n=1}^\infty \frac{z^n}{n^2} \qquad |z| \leq 1$$
which has integral form 
$$Li_2(z) = -\int_0^z \frac{\log(1-t)}{t}dt.$$

In \cite{Rog07}, Rogers introduced the following normalization. 
$$\Rog(z) = Li_2(z) + \frac{1}{2}\log |z| \log(1-z).$$
The dilogarithm function arises in many areas of mathematics, especially hyperbolic geometry and number theory (see \cite{Zag}).  In particular the volume of an ideal hyperbolic tetrahedron $T$ in $\Hs$ with ideal vertices $z_1,z_2,z_3,z_4 \in \hat\CC$ is 
$$Vol(T) = \Im(\Rog([z_1,z_2,z_3,z_4]))$$
where $[z_1,z_2,z_3,z_4]$ is the cross-ratio.

{\bf Pell's Equation:} Pell's equation for  $n \in \N$ is the Diophantine equation $x^2-ny^2 = \pm 1$ over $\Z$. Pell's equation has a long and interesting history going back to Archimedes' cattle problem (see  \cite{Lenstra}). The equation only has solutions for $n$ square-free, so we assume $n$ is square-free. Also by symmetry, we only need consider solutions with $x, y > 0$. A solution is positive/negative depending on whether $x^2-ny^2 =  1$ or $x^2-ny^2 = -1$. For all square-free $n$ there is always a positive solution but not necessarily a negative solution.  Solutions to Pell's equation correspond to units in $\Z[\sqrt{n}]$ by identifying $x,y$ with $x+y\sqrt{n}$ and  it is natural to identify the two. The smallest positive unit $u = x+y\sqrt{n}$ is called the {\em fundamental unit} and a well-known result gives that the set of positive units is exactly $\{u^k\}, k \in \N$ (see  \cite[Theorem 7.26]{MNZ})). 

In this paper we prove a new and surprising connection between the dilogarithm and solutions to Pell's equation. Using earlier work of the author which gave a dilogarithm identity associated to a hyperbolic surface, we obtain a dilogarithm identity for each solution $x,y$ to Pell's equation whose terms are given by the continued fraction expansion of $x+y\sqrt{n}$.

\subsection{Dilogarithm Identities}
The dilogarithm satisfies a number of classical identities. In particular by adding power series we have the squaring identity
$$Li_2(z) + L(-z) = \frac{1}{2}Li_2(z^2).$$ 
It follows by direct computation that this identity holds for the Rogers dilogarithm with
$$\Rog(z) +\Rog(-z) = \frac{1}{2}\Rog(z^2) \qquad\qquad\mbox{(Squaring Identity)}.$$
The other classic identities are the reflection identity of Euler
$$\Rog(x) + \Rog(1-x) = \frac{\pi^2}{6}  \qquad \Rog(x) + \Rog(x^{-1}) = \frac{\pi^2}{6} \qquad\qquad\mbox{(Reflection Identity)}$$
 Abel's well-known 5-term identity
$$\Rog(x) + \Rog(y) = \Rog(xy) + \Rog\left(\frac{x(1-y)}{1-xy}\right) + \Rog\left(\frac{y(1-x)}{1-xy}\right) \qquad\qquad\mbox{(Abel Identity)} $$
and Landau's identity
$$\Rog\left(-\frac{1}{x}\right) = -\Rog\left(\frac{1}{x+1}\right) \qquad  x > 0 \qquad\mbox{(Landau Identity)}.$$

A closed form for values of $\Rog$  is only known for a small set of values. These are
$$\Rog(0)  = 0 \qquad \Rog(1) = \frac{\pi^2}{6} \qquad \Rog\left(\frac{1}{2}\right) = \frac{\pi^2}{12}\qquad  \Rog(\phi^{-1}) = \frac{\pi^2}{10}\qquad \Rog(\phi^{-2}) = \frac{\pi^2}{15}$$
where $\phi$ is the golden ratio.  
In \cite{Lew91} Lewin gave the following remarkable infinite identity. 
 \begin{equation} \sum_{k=2}^\infty  \Rog\left(\frac{1}{k^2}\right) = \frac{\pi^2}{6}.
 \label{Lewin1}
 \end{equation}
  
 \section{Results}
 Using earlier work of the author we first prove the below new infinite identity for $\Rog$. We prove
\begin{theorem}
If $L > 0$ then
$$ \Rog(e^{-L}) = \sum_{k=2}^\infty  \Rog\left(\frac{\sinh^2\left(\frac L 2\right)}{\sinh^2\left(\frac{kL}{2}\right)}\right)$$
\label{Lid}
\end{theorem}
One immediate observation  is if we let $L \rightarrow 0$ we recover the formula of Lewin in equation \ref{Lewin1} above. 

We now apply the above identity to solutions of Pell's equation and units in the ring $\Z[\sqrt{n}]$.

\subsection*{Dilogarithm identity for Solution to Pell's Equation} In order to obtain our identity associated to a given solution $a^2-nb^2 = \pm 1$ of Pell's equation, we will choose $e^{L/2}  = a+b\sqrt{n}$ in the above. We then show that the righthandside is given in terms of the continued fraction expansion of $a+b\sqrt{n}$. We obtain:

  \begin{theorem}
Let $u = a+b\sqrt{n} \in \Z[\sqrt{n}]$ be a  solution to Pell's equation.
\begin{itemize}

\item If $u$ is a positive solution  with continued fraction convergents $r_j = h_j/k_j$ then    
$$\Rog\left(\frac{1}{u^2}\right) = \sum_{k=1}^\infty \Rog\left(\frac{1}{(h_{2k-1})^2} \right).$$
\item If $u$ is a negative solution and $u^2$ has convergents $R_j = H_j/K_j$ then 
$$\Rog\left(\frac{1}{u^2}\right) =  \sum_{k=0}^\infty \Rog\left(\frac{1}{b^2n(2H_{2k-1})^2} \right)+\Rog\left(\frac{1}{(2H_{2k+1}-H_{2k})^2} \right).$$
\end{itemize}
\label{fractions}\end{theorem}

\subsection*{Examples}
We now consider some examples:

{\bf Case of $\Z[\sqrt{2}]$:} 
For $\Z[\sqrt{2}]$ the fundamental unit is $3+2\sqrt{2}$ giving
$$\Rog\left(\frac{1}{(3+2\sqrt{2})^2}\right) = \Rog\left(\frac{1}{6^2} \right)+\Rog\left(\frac{1}{35^2} \right)+\Rog\left(\frac{1}{204^2} \right)+\Rog\left(\frac{1}{1189^2} \right)+\ldots.$$
We note that $3+2\sqrt{2}$ has convergents $r_k$ given by
$$\frac{5}{1},\frac{6}{1},\frac{29}{5}, \frac{35}{6}, \frac{169}{29},\frac{204}{35},\frac{985}{169},\frac{1189}{204}.$$

 It can be further shown that the units of $\Z[\sqrt{2}]$ are given by $(1+\sqrt{2})^k$.  As $u$ is a negative solution to Pell's equation and $u^2 =3+2\sqrt{2}$ we get
$$\Rog\left(\frac{1}{3+2\sqrt{2}}\right)=\Rog\left(\frac{1}{2(2)^2} \right)+\Rog\left(\frac{1}{7^2} \right)+\Rog\left(\frac{1}{2(12)^2} \right)+\Rog\left(\frac{1}{41^2} \right)+\Rog\left(\frac{1}{2(70)^2} \right)+\Rog\left(\frac{1}{239^2} \right)+\Rog\left(\frac{1}{2(408)^2} \right)+\ldots.$$

{\bf Case of $\Z[\sqrt{13}]$:} An interesting case of a large fundamental solution occurs for $\Z[\sqrt{13}]$. Here $u = 649+180\sqrt{13}$ is the fundamental unit giving
$$\Rog\left(\frac{1}{842401+233640\sqrt{13}}\right) = \Rog\left(\frac{1}{1298^2} \right)+\Rog\left(\frac{1}{1684803^2} \right)+\Rog\left(\frac{1}{2186872996^2}\right)\ldots.$$
The continued fraction convergents of $u$ are
$$\frac{1297}{1}, \frac{1298}{1},\frac{1683505}{1297},\frac{1684803}{1298},\frac{2185188193}{1683505},\frac{2186872996}{1684803}\ldots$$

\subsection*{Pell's equation over $\Q$} Similarly we  consider Pell's equation over $\Q$. If $a,b\in \Q$ satisfy Pell's equation $a^2-nb^2 =\pm 1$  we will identify this with the element $a+b\sqrt{n} \in \Q[\sqrt{n}]$. Applying the identity in Theorem \ref{Lid} we get the following.

\begin{theorem}
Let $u = a+ b\sqrt{n} \in \Q[\sqrt{n}], a, b > 0$ satisfy Pell's equation and let $u^k = a_k + b_k\sqrt{n}$.

 If $u$ is  a positive solution, then
$$\Rog\left(\frac{1}{u^2}\right) = \sum_{k=2}^\infty \Rog\left(\frac{1}{(b_k/b)^2} \right).$$
Further if $u \in \Z[\sqrt{n}]$ then $b_k/b \in \Z$ for all $k$.

If $u$ is a negative solution then  
$$\Rog\left(\frac{1}{u^2}\right) =  \sum_{k=1}^\infty \Rog\left(\frac{1}{n(b_{2k}/a)^2} \right)+\Rog\left(\frac{1}{(a_{2k+1}/a)^2} \right)$$
Further if $u \in \Z[\sqrt{n}]$ then $b_{2k}/a,a_{2k+1}/a \in \Z$ for all $k$.
\label{pell}
\end{theorem}

{\bf Fibonacci Numbers:}
The golden mean $\phi \in \Q[\sqrt{5}]$ corresponds to a negative solution to Pell's equation over $\Q$. Also we have
$$\phi^k = \frac{g_{k}+f_{k}\sqrt{5}}{2}$$
where $f_k$ is the classic Fibonacci sequence $1,1,2,3,5\ldots$ and $g_k$ is the Fibonacci sequence $1,3,4,7,11,\ldots$. 

As $\Rog(\phi^{-2}) = \pi^2/15$ we get the identity
$$\sum_{k=1}^\infty \left(\Rog\left(\frac{1}{5f_{2n}^2} \right)+\Rog\left(\frac{1}{g_{2n+1}^2} \right)\right) = \frac{\pi^2}{15}$$

\subsection*{Chebyshev Polynomials, Pell's Equation and Dilogarithms}
Chebyshev polynomials arise in numerous areas of mathematics and have a natural interpretation in terms of Pell's equation. The Chebyshev polynomials of the first kind $T_n$ are the unique polynomials satisfying $T_n(\cos(\theta)) = \cos(n\theta)$ and the Chebyshev polynomials of the second kind $U_n$ are given by
$$U_n(\cos(\theta)) = \frac{\sin((n+1)\theta)}{\sin(\theta)}$$
We obtain the following corollary. 
\begin{corollary}
Let $x > 1$ then
$$\Rog\left(\frac{1}{\left(x+\sqrt{x^2-1}\right)^2}\right) = \sum_{n=1}^\infty\Rog\left(\frac{1}{U_n(x)^2}\right).$$
\label{cheby}
\end{corollary}

    The  reader interested in knowing more about the dilogarithm function and its generalizations we refer  to the book \cite{Lew91}, {\em Structural Properties of Polylogarithms},  by L. Lewin and the aforementioned article \cite{Zag}, {\em The dilogarithm function},  by D. Zagier.

\section{Units in $\Z[\sqrt{n}]$, Pell's equation}
 We assume $n$ is not a perfect square. If $a+b\sqrt{n}\in \Z[\sqrt{n}]$ is a unit, then so are $\pm a\pm b\sqrt{n}$ and therefore we only need to consider  solutions $(a,b) \in \N^2$. 
It follows easily that $a\pm b\sqrt{n} \in\Z[\sqrt{n}]$  is a unit if and only if  $(a,b)$ satisfy  {\em Pell's equation} over $\Z$
$$a^2-nb^2=\pm 1.$$
We call a solution $(a,b)$ (or the unit $a+b\sqrt{n}$) positive/negative depending on if the righthandside of the Pell equation is positive/negative. Whereas there is always a solution to the positive Pell equation $x^2-ny^2 =1$, it can be shown that there are no solutions to $x^2-ny^2 = -1$ for certain $n$ (see \cite[Chapter 7]{MNZ}).

 \subsection*{Continued Fraction Convergents}
 
If  $u  \in \R_+$  we say $u$ has continued fraction expansion $u=[c_0,c_1,c_2,c_3,\ldots]$ if $c_i \in \Z$ and
$$u = c_0+\cfrac{1}{c_1+\cfrac{1}{c_2+\cfrac{1}{c_3 +\ldots}}}$$
 This means that if we define $r_n = [c_0,c_1,c_2,\ldots,c_n] \in \Q$ to be the $n^{th}$ {\em convergent}, then $r_n \rightarrow u$ as $n \rightarrow \infty$. If the continued fraction coefficients satisfy $c_{n+r} = c_{n}$ for $n > k$ we say $u$ is periodic with period $r$ and write $u =[c_0,c_1,\ldots,c_{k},\overline{c_{k+1},\ldots,c_{k+r}}]$. We have the following  standard description of $r_n$;
 \begin{theorem}{(\cite[Theorems 7.4, 7.5]{MNZ})}
 Let $u \in \R_+$ with $u = [c_0,c_1,c_2,\ldots]$ and define $h_n,k_n$ by 
 $$h_i = c_ih_{i-1}+h_{i-2}\qquad k_i = c_ik_{i-1}+k_{i-2}\qquad  i \geq 0$$
 with $(h_{-2},k_{-2}) = (0,1), (h_{-1}, k_{-1}) = (1,0)$. Then $\gcd(h_i,k_i) = 1$ and
 $$r_n = [c_0,c_1,c_2,\ldots, c_n] = \frac{h_n}{k_n}$$
 \end{theorem}

  The positive units in $\Z[\sqrt{n}]$ have the following elegant description.
 \begin{theorem}{(\cite[Theorem 7.26]{MNZ})}
Let $n \in \N$ not be a perfect square. Then there is a unique solution $(a,b) \in \N^2$  of Pell's equation $x^2-ny^2 =1$ such that the set of solutions to $x^2-ny^2 =1$ in $\N^2$ is $\{(a_k,b_k)\}_{k=1}^\infty$ where
$$a_k + b_k\sqrt{n} = (a+b\sqrt{n})^k.$$
\label{posunits}
\end{theorem}

The pair $(a,b)$ is called the {\em fundamental solution}  of $x^2-ny^2 = 1$.  Thus one consequence of the above is if we let  $u $ be the fundamental unit then $\{u^k\}$ gives the set of all positive solutions to Pell's equation and the dilogarithm identity in Theorem \ref{fractions} can be interpreted as a sum over all solutions to Pell's equation.

\section{The Orthospectrum Identity}
In a prior paper, the author proved a dilogarithm identity for a hyperbolic surface with geodesic boundary. The identity was generalized to hyperbolic manifolds by the author and Kahn in  \cite{BK10}.
The relation to other identities on hyperbolic manifolds such as the Basmajian identity (see \cite{Bas93}), McShane-Mirzakhani identity (see \cite{McS98}, \cite{Mir07}) and Luo-Tan identity (see \cite{LT11}) is discussed in \cite{BT16}.

In order to state the orthospectrum identity, we recall  some basic terms.

If $S$ is a hyperbolic surface with totally geodesic boundary, an {\em orthogeodesic} $\alpha$ is a proper geodesic arc which is perpendicular to the boundary $\partial S$ at its endpoints. The set of orthogeodesics of $S$ is denoted $O(S)$. Each boundary component is either a closed geodesic or an infinite geodesic whose endpoints are {\em boundary cusps} of $S$. We let $N(S)$ be the number of boundary cusps of $S$. 
Further let $\chi(S)$ be given by $Area(S) = 2\pi|\chi(S)|$.

One elementary example of a surface will be an ideal n-gon which has $N(S) = n$ and  $O(S)$ a finite set and $\chi(S) = 1-n/2$. In fact these are the only surfaces with $O(S)$ finite.

The dilogarithm orthospectrum identity is as follows;

\begin{theorem}{(Dilogarithm Orthospectrum Identity, \cite{B11})}
Let $S$ be a finite area hyperbolic surface with totally geodesic boundary $\partial S \neq 0$. Then
$$\sum_{\alpha \in O(S)} \Rog\left(\frac{1}{\cosh^2\left(\frac{l(\alpha)}{2}\right)}\right) = -\frac{\pi^2}{12}(6\chi(S) + N(S))$$ 
and equivalently
$$\sum_{\alpha \in O(S)} \Rog\left(-\frac{1}{\sinh^2\left(\frac{l(\alpha)}{2}\right)}\right) = \frac{\pi^2}{12}(6\chi(S) + N(S))$$ 
\label{id1}
\end{theorem}

In the original paper \cite{B11}, we showed that the above identity  recovers the reflection identity, Abels identity and Landau's identity by considering the elementary cases of the ideal quadrilateral and ideal pentagon respectively.

\section{An infinite dilogarithm identity}
We define the cross-ratio of 4 distinct points in $\hat\CC$ by
$$[z_1,z_2,z_3,z_4] =\frac{(z_1-z_2)(z_4-z_3)}{(z_1-z_3)(z_4-z_2)}.$$
We let $x_1,x_2,x_3,x_4 \in \partial \Hp$ be distinct points ordered counterclockwise on $\partial \Hp$. If  $g$ is the geodesic with endpoints $x_1,x_2$ and $h$ is the geodesic with endpoints $x_3,x_4$ then $g,h$ are disjoint. A simple calculation shows that the perpendicular distance between them $l$ is given by
$$\frac{1}{\cosh^2(l/2)} = [x_1,x_2,x_3,x_4].$$

\begin{figure}[htbp] 
   \centering
   \hspace{1.5in}\includegraphics[width=4in]{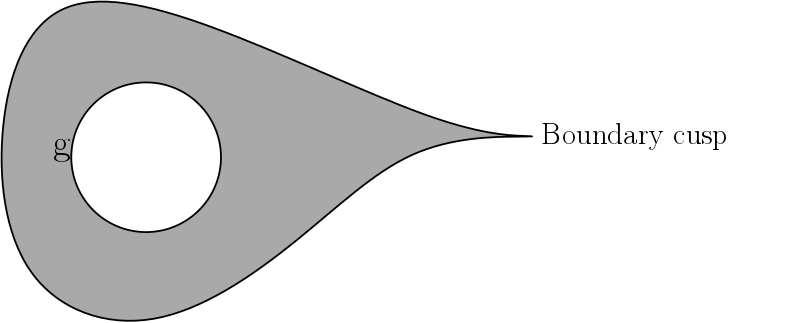} 
   \caption{Surface S}
   \label{fig1}
\end{figure}
We now  prove Theorem \ref{Lid} which we now restate.
\vspace{.25in}

\noindent{\bf Theorem \ref{Lid}}
{\em If $L > 0$ then
$$ \Rog(e^{-L}) = \sum_{k=2}^\infty  \Rog\left(\frac{\sinh^2\left(\frac L 2\right)}{\sinh^2\left(\frac{kL}{2}\right)}\right)$$
}
{\bf Proof:}
We consider the hyperbolic surface $S$ which is topologically an annulus with one boundary component being a closed geodesic $g$ of length $L$ and the other $h$ an infinite geodesic with a single boundary cusp (see figure \ref{fig1}).

\begin{figure}[htbp] 
    \centering
    \includegraphics[width=4in]{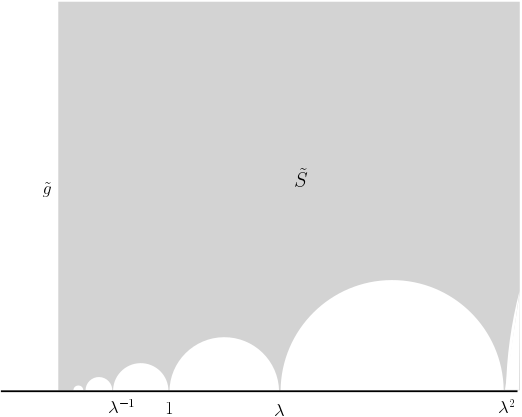} 
    \caption{Universal cover of $S$}
    \label{fig2}
 \end{figure}

 We lift $S$ to the upper half plane with $g$ lifted to the y-axis and let $\lambda = e^L$. Then $\tilde S$ is an infinite sided ideal polygon, which is invariant under multiplication by $\lambda$ (see figure \ref{fig2}). 
 
 We  normalize so that one of the ideal vertices is at $1$. Then the vertices of $\tilde S$ are $\lambda^k$ for $k \in \Z$. There is a single orthogeodesic with endpoint on $g$ and it has length $l$ satisfying
$$\frac{1}{\cosh^2(l/2)} = [\infty,0,1,\lambda] = \frac{\lambda-1}{\lambda} = 1-e^{-L}.$$
The other ortholengths are given by $l_k$ where
$$\frac{1}{\cosh^2(l_k/2)} = [1,\lambda, \lambda^k,\lambda^{k+1}] = \frac{(1-\lambda)(\lambda^{k+1}-\lambda^k)}{(1-\lambda^k)(\lambda^{k+1}-\lambda)} =\lambda^{k-1}\frac{(\lambda-1)^2}{(\lambda^k-1)^2}$$
$$=\frac{(\lambda^{1/2}-\lambda^{-1/2})^2}{(\lambda^{k/2}-\lambda^{-k/2})^2}= \frac{\sinh^2(L/2)}{\sinh^2(kL/2)}.$$
Thus the dilogarithm identity for $S$ is
$$\Rog(1-e^{-L}) +\sum_{k=2}^\infty \Rog\left(\frac{\sinh^2(L/2)}{\sinh^2(kL/2)}\right) = \frac{\pi^2}{6}$$
Using the reflection identity $\Rog(1-x) + \Rog(x) = \pi^2/6$ we get
$$\Rog(e^{-L}) = \sum_{k=2}^\infty \Rog\left(\frac{\sinh^2(L/2)}{\sinh^2(kL/2)}\right).$$
\eproof

 \section{Proof of identity for Solutions to Pell's Equation over $\Q$}
We now prove the dilogarithm identity for solutions to Pell's equation over $\Q$ given in Theorem \ref{pell}. 

\vspace{.15in}

{\bf Proof of Theorem \ref{pell}:}

Let   $e^{L/2} = u = a+b\sqrt{n}$, then $e^{-L/2}= u^{-1} =\pm(a-b\sqrt{n})$ with the sign depending on if $u$ is a positive unit or negative. If $u$ is a positive unit then
$$\cosh(L/2) =  a \qquad \sinh(L/2) = b\sqrt{n} $$
and if $u$ is a negative unit
$$\sinh(L/2) =  a \qquad \cosh(L/2) = b\sqrt{n} $$
Either way we have 
$$u^k  = e^{kL/2} = \cosh(kL/2) + \sinh(kL/2).$$
We let $m_k =\sinh(kL/2)$ and $n_k = \cosh(kL/2)$.  
The dilogarithm identity gives 
$$\Rog\left(\frac{1}{u^2}\right) = \sum_{k=2}^\infty \Rog\left( \frac{\sinh^2(L/2)}{\sinh^2(kL/2)}\right) = \sum_{k=2}^\infty \Rog\left( \frac{m^2_1}{m_k^2}\right).$$
If $u$ is a positive root,  then $m_1 = \sinh(L/2) = b\sqrt{n}$ and $n_1 = \cosh(L/2) =  a$.
Then by the addition formulae we have
$$m_{k+1} =a.m_k + bn_k\sqrt{n}  \qquad n_{k+1} = n_ka +bm_k\sqrt{n}.$$
Then by induction we have $n_k = a_k$ and $m_k = b_k\sqrt{n}$ and 
$$b_{k+1} = ab_k+ba_k \qquad a_{k+1} = aa_k+nbb_k.$$
$$\Rog\left(\frac{1}{u^2}\right) = \sum_{k=1}^\infty \Rog\left(\frac{b^2}{b_k^2} \right) = \sum_{k=1}^\infty \Rog\left(\frac{1}{(b_k/b)^2} \right).$$

If $u$ is a negative solution,  then $m_1 = \sinh(L/2) = a$ and $n_1 = \cosh(L/2) =  b\sqrt{n}$.
Then by the addition formulae we have
$$m_{k+1} = bm_k\sqrt{n} +an_k \qquad n_{k+1} = bn_k\sqrt{n} +am_k.$$
Therefore 
$$n_{2k} = a_{2k} \qquad n_{2k+1} = b_{2k+1}\sqrt{n} \qquad m_{2k} = b_{2k}\sqrt{n} \qquad m_{2k+1} = a_{2k+1}.$$
Therefore
$$b_{2k} = ba_{2k-1}+ab_{2k-1}\qquad a_{2k+1} = bb_{2k}n+aa_{2n}.$$
Therefore  
$$\Rog\left(\frac{1}{u^2}\right) = \sum_{k=1}^\infty \Rog\left(\frac{m_1^2}{m_k^2} \right) = \sum_{k=1}^\infty \Rog\left(\frac{1}{n(b_{2k}/a)^2} \right)+\Rog\left(\frac{1}{(a_{2k+1}/a)^2} \right).$$
 \eproof

We now prove Corollary \ref{cheby} relating the identity to the Chebyshev polynomials $U_n$ of the second kind.

{\bf Proof of Corollary \ref{cheby}:}
We have the Chebyshev polynomials $T_n(x), U_n(x) \in\R[x]$. We let $x = \cos(\theta)$, then $\sin(\theta) = \sqrt{1-x^2}$. Therefore
$$e^{i\theta} = \cos(\theta) + i\sin(\theta) = x+i\sqrt{1-x^2} = x+\sqrt{x^2-1}$$
and 
$$e^{in\theta} = (x+\sqrt{x^2-1})^n\qquad e^{-in\theta} = (x-\sqrt{x^2-1})^n.$$
 Thus
 $$T_n(x) = \cos(n\theta) = \frac{1}{2}\left((x+\sqrt{x^2-1})^n+(x-\sqrt{x^2-1})^n\right)$$
and
$$U_{n-1}(x) = \frac{\sin(n\theta)}{\sin{\theta}} = \frac{1}{2\sqrt{x^2-1}} \left((x+\sqrt{x^2-1})^n-(x-\sqrt{x^2-1})^n\right).$$
As this holds for $|x| < 1$, it also holds for  all $x \in \R$.
Now if $x > 1$ then we let $x = \cosh(L/2)$  then $\sqrt{x^2-1} = \sinh(L/2)$ giving
$$x+\sqrt{x^2-1} = e^{L/2} \qquad x-\sqrt{x^2-1} = e^{-L/2}.$$
Therefore by the above formulae
$$T_k(x) = \frac{e^{kL/2}+e^{-kL/2}}{2} = \cosh(kL/2) \qquad U_{k-1}(x) =\frac{e^{kL/2}-e^{-kL/2}}{2\sinh(L/2)} = \frac{\sinh(kL/2)}{\sinh(L/2)}.$$
Thus
$$\Rog\left(\frac{1}{\left(x+ \sqrt{x^2-1}\right)^2}\right) = \sum_{k=2}^\infty\Rog\left(\frac{\sinh^2(L/2)}{\sinh^2(kL/2)}\right) =  \sum_{k=1}^\infty\Rog\left(\frac{1}{U_{k}(x)^2}\right).$$
\eproof

 \section{Identity for continued fraction convergents}
We now  consider the case where $u \in \Z[\sqrt{n}]$ and prove Theorem \ref{fractions} expressing the above in terms of the convergents $r_j = h_j/k_j$ of their continued fractions expansion. First we have the following lemma.

\begin{lemma}
Let $u =a + b\sqrt{n} \in \Z[\sqrt{n}]$ be a solution to Pell's equation with $a,b \in \N$.
If $u$ is a positive solution then $u =  [2a-1,\overline{1,2a-2}]$ and if $u$ is a negative solution then $u = [\overline{2a}]$.
\end{lemma}
{\bf Proof:}
if $u$ is a negative solution. then $u = a + \sqrt{a^2+1}$ and $u^2-2au-1 = 0$.
Then
$$u = 2a+\frac{1}{u}$$
 giving $u =  [\overline{2a}]$.
 
 If $u$ is a positive solution then $u = a + \sqrt{a^2-1}$. Therefore $u$ satisfies the quadratic
$u^2 -2au+1 = 0$. Rewriting we have
$$u = 2a-\frac{1}{u} = 2a-1 + 1-\frac{1}{u} = 2a-1 +\frac{u-1}{u}$$
Now we have
$$\frac{u-1}{u} = \frac{1}{\frac{u}{u-1}} = \frac{1}{1 + \frac{1}{u-1}} = \frac{1}{1 + \frac{1}{2a-2 +\frac{u-1}{u}}}.$$
Therefore $u =  [2a-1,\overline{1,2a-2}]$. 
 
\eproof

Using the above description of the continued fraction, we will show the relation between the approximates $r_j = h_j/k_j$ for $u$ and the coefficients $a_j,b_j$ given by $u^j = a_j +b_j\sqrt{n}$. This will allow us to prove Theorem \ref{fractions} giving the dilogarithm identity in terms of the convergents of the the continued fraction convergents.

\begin{lemma}
Let $u = a+b\sqrt{n} \in \Z[\sqrt{n}]$ be a  solution to Pell's equation.  

If $u$ is a positive solution and $u$ has continued fraction convergents $r_j = h_j/k_j$ then  $k_{j} = h_{j-2}$ and 
$$\Rog\left(\frac{1}{u^2}\right)  = \sum_{j=1}^\infty \Rog\left(\frac{1}{(h_{2j-1})^2}\right) .$$
If $u$ is a negative solution and $u^2$ has continued fraction convergents $R_j = H_j/K_j$ then 
$$\Rog\left(\frac{1}{u^2}\right)  = \sum_{j=0}^\infty \left(\Rog\left(\frac{1}{nb^2(2H_{2j-1})^2}\right)+\Rog\left(\frac{1}{(2H_{2k+1} - H_{2k})^2}\right)\right) .$$

\end{lemma}

{\bf Proof:}
Let $u = a+b\sqrt{n} = e^{L/2}$, then $u^k = a_k + b_k\sqrt{n} = \cosh(kL/2)+\sinh(kL/2)$.
As $u$ is a positive solution  $u = [2a-1,\overline{1,2a-2}]$. Therefore  we have $(h_0,h_{-1}) = (2a-1,1)$ and for $k > 0$
$$\begin{bmatrix} 
h_{2k}\\
h_{2k-1}
\end{bmatrix}
=
A^{k}
\begin{bmatrix} 
2a-1\\
1
\end{bmatrix}
\qquad\mbox{where}\qquad A = 
 \begin{bmatrix} 
2a-2 & 1\\
1 & 0
\end{bmatrix}
 \begin{bmatrix} 
1 & 1\\
1 & 0
\end{bmatrix}
= \begin{bmatrix} 
2a-1 & 2a-2\\
1 & 1
\end{bmatrix}
.$$
The matrix $A$ has characteristic polynomial $x^2-2ax+1$ giving eigenvalues $u, 1/u$ and eigenvectors $(u-1,1), (1-u, u)$. Therefore diagonalizing we get
$$\begin{bmatrix} 
h_{2k}\\
h_{2k-1}
\end{bmatrix}
=\frac{1}{u^2-1}\begin{bmatrix} 
u-1 & 1-u\\
1 & u
\end{bmatrix}
\begin{bmatrix} 
u^{k} & 0\\
0 & u^{-k}
\end{bmatrix}
\begin{bmatrix} 
u & u-1\\
-1 & u-1
\end{bmatrix}
\begin{bmatrix} 
2a-1\\
1
\end{bmatrix}
.$$
Multiplying out and noting that $u = e^{L/2}$ we have
\begin{eqnarray}
h_{2k} &=& \frac{(u-1)\left(u^{k+2}+u^{-(k+1)}\right)}{u^2-1} = \frac{\cosh((k+\frac{3}{2})L/2)}{\cosh(L/4)}\label{heqn1}\\
h_{2k-1} &=& \frac{u^{k+2}-u^{-k}}{u^2-1} =  \frac{\sinh((k+1)L/2)}{\sinh(L/2)}.
\label{heqn2}
\end{eqnarray}

It follows that  for $k\geq 1$
\begin{equation}
h_{2k-3}  = \frac{\sinh(kL/2)}{\sinh(L/2)} = \frac{b_k}{b} .
\label{eq_pos}
\end{equation}

Therefore
$$\Rog\left(\frac{1}{u^2}\right) = \sum_{k=2}^\infty \Rog\left(\frac{1}{(b_k/b)^2}\right) = \sum_{j=1}^\infty \Rog\left(\frac{1}{(h_{2j-1})^2}\right) .$$
Similarly we note that as $(k_0, k_1) = (1,0)$ then applying the above analysis we get
$$k_{2j} = \frac{\cosh((j+\frac{1}{2})L/2)}{\cosh(L/4)}  = h_{2j-2}$$
and 
$$k_{2j-1} = \frac{\sinh(jL/2)}{\sinh(L/2)} = h_{2j-3}$$
giving $k_{j} = h_{j-2}$.

If $u$ is a negative solution, then  for $k$ odd $a_{k} =  \sinh(kL/2), b_{k}\sqrt{n} = \cosh(kL/2)$ and  for $k$ even
$b_{k}\sqrt{n} =  \sinh(kL/2), a_{k} = \cosh(kL/2)$.

 As $u = [\overline{2a}]$ we have the formula
$$h_{j+1} = 2a h_{j}+h_{j-1} \qquad k_{j+1} = 2ak_j+k_{j+1}$$ 
with $(h_{-2},k_{-2}) = (0,1)$ and $(h_{-1},k_{-1}) = (1,0)$.
Iterating we get $h_j = 0,1,2a,\ldots$ and $k_j = 1,0, 1, 2a, \ldots$. Therefore $k_j = h_{j-1}$ for $j \geq -1$. We focus on calculating $h_k$. As $(h_{-1},h_{-2}) = (1,0)$ we have the recursion
$$\begin{bmatrix} 
h_k\\
h_{k-1}
\end{bmatrix}
=
A^{k+1}
\begin{bmatrix} 
1\\
0
\end{bmatrix}
\qquad \mbox{where } \qquad A = \begin{bmatrix} 
2a & 1\\
1 & 0
\end{bmatrix}.$$
The matrix $A$ has characteristic polynomial $x^2-2ax-1 =  0$ giving eigenvalues $u, -u^{-1}$ and eigenvectors $(u,1), (1,-u)$. Thus
$$\begin{bmatrix} 
h_k\\
h_{k-1}
\end{bmatrix}
=\frac{1}{u^2+1}\begin{bmatrix} 
u & 1\\
1 & -u
\end{bmatrix}
\begin{bmatrix} 
u^{k+1} & 0\\
0 & (-u)^{-k-1}
\end{bmatrix}
\begin{bmatrix} 
u & 1\\
1 & -u
\end{bmatrix}
\begin{bmatrix} 
1\\
0
\end{bmatrix}
.$$
Multiplying we get
$$ h_k = \frac{1}{u^2+1}\left( u^{k+3}+ (-1)^{k+1} u^{-(k+1)}\right) = \frac{1}{u+u^{-1}}\left( u^{k+2}+ (-1)^{k+1} u^{-(k+2)}\right).$$
For $k$ odd we have
$$h_k = \frac{\cosh((k+2)L/2)}{\cosh(L/2)} = \frac{b_{k+2}}{b}.$$
Similarly for $k$ even we have
$$h_k = \frac{\sinh((k+2)L/2)}{\cosh(L/2)} = \frac{b_{k+2}}{b}.$$
Thus for all $k \geq 0$
$$\frac{b_k}{b} = h_{k-2}.$$
We let $H_j,K_j$ be the convergents for the continued fraction expansion of $u^2$.  Then as $u^2 = e^{L}$ is a positive solution to Pell's equation. Applying  equations \ref{heqn1}  and \ref{heqn2} above we have,
$$H_{2k} = \frac{\cosh((2k+3)L/2}{\cosh(L/2)} = \frac{b_{2k+3}}{b} = h_{2k+1}.$$
$$H_{2k-1} = \frac{\sinh((k+1)L)}{\sinh(L)} = \frac{1}{2}\left(\frac{\sinh(2k+1)L/2}{\sinh(L/2)} + \frac{\cosh(2k+1)L/2}{\cosh(L/2)} \right) = \frac{1}{2}\left(\frac{a_{2k+1}}{a}+h_{2k-1} \right).$$
Also if $(u^2)^k = A_k+B_k\sqrt{n}$ then $A_k = a_{2k}, B_{k} = b_{2k}$ and by equation \ref{eq_pos}
$$H_{2k-3} = \frac{B_{k}}{B_1} = \frac{b_{2k}}{2ab} = \frac{h_{2k-2}}{2a}$$
giving
$$h_{2k} = 2aH_{2k-1}\qquad h_{2k+1} = H_{2k}.$$
Also
$$\frac{b_{2k}}{a} = 2bH_{2k-3}\qquad \frac{a_{2k+1}}{a} = 2H_{2k-1}-h_{2k-1} = 2H_{2k-1} - H_{2k-2}.$$
Thus if $u$ is a negative solution to Pell's equation
$$\Rog\left(\frac{1}{u^2}\right) = \sum_{k=0}^\infty \Rog\left(\frac{1}{n(2bH_{2k-1})^2}\right)+ \Rog\left(\frac{1}{(2H_{2k-1} - H_{2k-2})^2}\right).$$
\eproof

\section{Ideal n-gon identities}

 Ramanujan gave a number of {\em value-identities} for linear combinations of specific values of $\Rog$ as follows (see \cite[Entry 39]{Ram1}); 
\begin{enumerate}
\item 
$$Li_2\left(\frac{1}{3}\right) - \frac{1}{6}Li_2\left(\frac{1}{9}\right)= \frac{\pi^2}{18}-\frac{\log^2 3}{6}$$
\item$$Li_2\left(-\frac{1}{2}\right)+\frac{1}{6}Li_2\left(\frac{1}{9}\right) = -\frac{\pi^2}{18} + \log 2\log 3-\frac{\log^2 2}{2} -\frac{\log^2 3}{3}$$
\item $$Li_2\left(\frac{1}{4}\right)+\frac{1}{3}Li_2\left(\frac{1}{9}\right) = \frac{\pi^2}{18} + 2\log 2\log 3-2\log^2 2 -\frac{2\log^2 3}{3}$$
\item $$Li_2\left(-\frac{1}{3}\right) - \frac{1}{3}Li_2\left(\frac{1}{9}\right)= -\frac{\pi^2}{18}-\frac{\log^2 3}{6}$$
\item $$Li_2\left(-\frac{1}{8}\right) + Li_2\left(\frac{1}{9}\right)=  -\frac{\log^2 9/8}{2}$$
\end{enumerate}
More recently in their article \cite{64id}, Bailey, Borwein, and Plouffe  gave the identity
\begin{equation}
36Li_2\left(\frac{1}{2}\right) -36Li_2\left(\frac{1}{4}\right)- 12Li_2\left(\frac{1}{8}\right) + 6Li_2\left(\frac{1}{64}\right) = \pi^2.
\label{64}
\end{equation}

Applying Landau's identity gives $\Rog(-1/3) = -\Rog(1/4)$ and $\Rog(-1/8) = -\Rog(1/9)$, which reduces the value-identities of Ramanujan  to the two equations
$$\Rog\left(\frac{1}{4}\right) + \frac{1}{3}\Rog\left(\frac{1}{9}\right) = \frac{\pi^2}{18} \qquad \Rog\left(\frac 1 3\right) - \frac{1}{6}\Rog\left(\frac 1 9\right) = \frac{\pi^2}{18}.$$

In this section we will show that these value-identities follow from considering ideal hyperbolic hexagons.

We recall the dilogarithm identity in \cite{B11} for ideal hyperbolic polygons. Let  $P$ be an ideal polygon in $\Hp$ with vertices in counterclockwise order $x_1,\ldots,x_n$ about $\partial \Hp$. If $l_{ij}$ is the length of the orthogeodesic joining side $[x_i,x_{i+1}]$ to $[x_j,x_{j+1}]$ then a simple calculation gives
$$[x_i,x_{i+1},x_j,x_{j+1}] = \frac{1}{\cosh^2(l_{ij}/2)}.$$

Applying the orthospectrum identity in Theorem \ref{id1} to $P$ we obtain the equation
$$\sum_{|i-j| \geq 2} \Rog([x_i,x_{i+1},x_j,x_{j+1}]) = \frac{(n-3)\pi^2}{6}.$$

If $P$ is the regular ideal n-gon then we obtain the equation
$$\frac{e_n}{2}\Rog\left(\sin^2\left(\pi/n\right)\right) +\sum_{k=2}^{\lfloor\frac{n}{2}\rfloor} \Rog\left(\frac{\sin^2(\pi/n)}{\sin^2(k\pi/n)}\right) = \frac{(n-3)\pi^2}{6n}$$
where $e_n = 0$ if n is odd and $e_n = 1$ if n is even.

\section{Ideal Hexagons and Ramanujan's value identities}
We now show that Ramanujan's value identities 1-5 and  identity \ref{64} of Bailey, Borwein, Plouffe, correspond to identities for the regular ideal hexagons.

For the regular 6-gon $H_{reg}$ the orthospectrum identity gives 
$$6\Rog\left(\frac 1 3\right) + 3\Rog\left(\frac 1 4\right) = \frac{\pi^2}{2}.$$
We also have from Landau's identity that $\Rog(-1/3) = -\Rog(1/4)$. Therefore applying the squaring identity we get
$$\frac{1}{2} \Rog\left(\frac 1 9\right) = \Rog\left(\frac 1 3\right) + \Rog\left(-\frac 1 3\right) = \Rog\left(\frac 1 3\right) - \Rog\left(\frac 1 4\right).$$
Thus we obtain
$$\Rog\left(\frac 1 3\right) - \Rog\left(\frac 1 4\right) = \frac{1}{2} \Rog\left(\frac 1 9\right).$$
Combining this and the identity above for the regular hexagon, we obtain Ramanujan's value-identities
$$\Rog\left(\frac{1}{4}\right) + \frac{1}{3}\Rog\left(\frac{1}{9}\right) = \frac{\pi^2}{18} \qquad \Rog\left(\frac 1 3\right) - \frac{1}{6}\Rog\left(\frac 1 9\right) = \frac{\pi^2}{18}.$$

To recover the  identity \ref{64} of Bailey, Borwein, Plouffe  we note that by Landau $\Rog(-1/8) = -\Rog(1/9)$  then by the squaring identity we have   
$$\frac{1}{2}\Rog\left(\frac{1}{64}\right) = \Rog\left(\frac{1}{8}\right) +\Rog\left(-\frac{1}{8}\right) = \Rog\left(\frac{1}{8}\right) -\Rog\left(\frac{1}{9}\right)$$
Therefore  substituting for $\Rog(1/8)$ we get
$$36\Rog\left(\frac{1}{2}\right) -36\Rog\left(\frac{1}{4}\right)- 12\Rog\left(\frac{1}{8}\right) + 6\Rog\left(\frac{1}{64}\right) = 36\Rog\left(\frac{1}{2}\right) -36\Rog\left(\frac{1}{4}\right)- 12\Rog\left(\frac{1}{9}\right).$$
As $\Rog(1/2) = \pi^2/12$ and by the hexagon identity $3\Rog(1/4)+\Rog(1/9) = \pi^2/6$ we 
recover identity \ref{64}.
 
\end{document}